\documentclass[12pt]{article}
\usepackage{amssymb}
\newcommand{\R}{\mathbb{R}}

\newcommand{\T}{\mathbb{T}}
\def\build#1_#2^#3{\mathrel{
\mathop{\kern 0pt#1}\limits_{#2}^{#3}}}
\def\llbracket{[\hspace{-.10em} [ }
\def\rrbracket{ ] \hspace{-.10em}]}
\def\cq{$\hfill \square$}

\def\cc{{\cal C}}

\def\m{{\cal M}}

\def\t{{\cal T}}

\def\m{{\bf m}}

\def\be{\begin{equation}}
\def\ee{\end{equation}}
\def\ba{\begin{eqnarray*}}
\def\ea{\end{eqnarray*}}
\def\ov{\overline}

\def\wt{\widetilde}

\def\la{\longrightarrow}
\def\da{\downarrow}
\def\ua{\uparrow}

\def\proof{\vskip 3mm \noindent{\bf Proof:}\hskip10pt}

\newtheorem{theorem}{Theorem}[section]
\newtheorem{lemma}[theorem]{Lemma}
\newtheorem{proposition}[theorem]{Proposition}
\newtheorem{corollary}[theorem]{Corollary}

\begin{document}

\title{ \bf The Hausdorff measure of stable trees}
\author{Thomas {\sc Duquesne} and
Jean-Fran\c cois {\sc Le Gall}}
\vspace{2mm}
\date{\tiny\today}

\maketitle

\begin{abstract} We study fine properties of the so-called
stable trees, which are the scaling limits of critical Galton-Watson trees
conditioned to be large. In particular we derive the exact Hausdorff
measure function for Aldous' continuum random tree and for its level sets.
It follows that both the uniform measure on the tree and the local
time measure on a level set coincide with certain Hausdorff measures.
Slightly less precise results are obtained for the Hausdorff measure
of general stable trees.
\end{abstract}

\section{Introduction}

The purpose of the present work is to study the Hausdorff measure
properties of the continuous random trees called stable trees.
Roughly speaking, stable trees are the continuous random trees that
arise as scaling limits of Galton-Watson trees with critical offspring
distribution, which are
conditioned to be large in some sense. In the most important case
where the offspring distribution also has a finite variance, this leads
to Aldous' continuum random tree (the CRT, see \cite{Al1} and
\cite{Al2}) and variants of the
CRT.  Alternatively, stable trees can be viewed as describing the
genealogical structure
of continuous-state branching processes with a stable branching mechanism
of the type $\psi(u)=u^\alpha$ for $1<\alpha\leq 2$. Thus they also
encode the genealogy of superprocesses with stable branching mechanism,
which have been studied by many authors. The case $\alpha=2$ yields
the so-called quadratic branching mechanism, corresponding to
finite variance superprocesses.

Stable trees are particular instances of the more general L\'evy trees
studied in \cite{DuLG}. In the formalism of \cite{DuLG}, L\'evy trees are
random variables taking values in the space of all (compact) rooted $\R$-trees.
Informally
an $\R$-tree is a metric space $(\t,d)$ such that for any two points
$\sigma$ and $\sigma'$
in $\t$ there is a unique arc with endpoints $\sigma$ and $\sigma'$
and furthermore this
arc is isometric to a compact interval of the real line. A rooted $\R$-tree
is an $\R$-tree with a distinguished vertex $\rho$ called the root.
We write $H(\t)$
for the height of $\t$, that is the maximal distance from the root to
a vertex in $\t$.
Two rooted $\R$-trees are called equivalent if there is a
root-preserving isometry
that maps one onto the other. It was noted in
\cite{EPW} that the set $\T$ of equivalence classes of compact rooted
$\R$-trees, equipped
with the Gromov-Hausdorff distance \cite{Gro}, is a Polish space.

It is shown in \cite{DuLG} that with every critical or subcritical
branching mechanism
function $\psi$ such that the corresponding branching process dies
out a.s. one can
associate a $\sigma$-finite measure $\Theta$ on $\T$ which is called the
``law'' of the L\'evy tree with branching mechanism $\psi$. Although $\Theta$
is an infinite measure, the quantity
$$v(\varepsilon)=\Theta(H(\t)>\varepsilon)$$
is finite for every $\varepsilon>0$ and is determined by the equation
$\int_{v(\varepsilon)}^\infty \psi(u)^{-1}du=\varepsilon$.
L\'evy trees enjoy
the important ``branching property'', which is analogous to a
classical result for
Galton-Watson trees: For every $a>0$, under
the probability measure $\Theta(\cdot\mid H(\t)>a)$ and conditionally
given the part of the
tree below level
$a$, the subtrees above that level are distributed as the atoms of a
Poisson point measure
whose intensity is a random multiple of $\Theta$ (the random factor is
the total mass of the local time measure at level $a$ that will be
discussed below).
It has recently been shown by Weill \cite{W} that this branching property
characterizes L\'evy trees.

When
$\psi(u)=u^\alpha$ for some $\alpha\in(1,2]$ we write
$\Theta_\alpha=\Theta$ and call
$\Theta_\alpha$ the law of the stable tree with index $\alpha$. In
addition to the branching
property, stable trees  possess the following scaling property. For
every $r>0$ and every tree
$\t\in\T$, denote by $r\t$ the ``same'' tree $\t$ with metric $d$
replaced by $r\,d$.
Then, for every $r>0$, the law of $r\t$ under $\Theta_\alpha(d\t)$ is
$r^{1\over
\alpha-1}\Theta_\alpha$.

An explicit construction
of $\Theta_\alpha$ may be given through the coding of real trees from the
height process studied in \cite{LGLJ1} and \cite{DuLG0} (see also Theorem 2.1
of \cite{DuLG} for the coding of real trees). This construction is especially
simple in the case $\alpha=2$, since the height process is then just a
Brownian excursion, and this approach essentially reduces to
Aldous' construction of the CRT from the normalized Brownian excursion
(Corollary 22 in \cite{Al2}). Alternatively, we may use the following
approximation
by discrete trees. Let $\pi$
be a probability distribution on $\{0,1,\ldots\}$. Assume that $\pi$
has mean $1$
and is in the domain of attraction of a stable distribution with
index $\alpha$,
in the sense that there exists an increasing sequence $(a_n)_{n=1,2,\ldots}$
of positive integers such that, if $\xi_1,\xi_2,\ldots$ are i.i.d. with
distribution $\pi$,
$(a_n)^{-1}(\xi_1+\cdots+\xi_n-n)$ converges  in distribution to a
stable distribution with
index $\alpha$. Let $c>0$ be a constant and for every $n\geq 1$ let $\theta_n$
be a Galton-Watson tree with offspring distribution $\pi$ conditioned
to have height
greater than $cn$. Notice that $\theta_n$ can be viewed as a random
$\R$-tree by
affecting length $1$ to each edge. Then the distribution of $n^{-1}\theta_n$
converges as $n\to \infty$ to the probability measure
$\Theta_\alpha(\cdot\mid H(\t)>c)$.
This result follows from a special case of Proposition 2.5.2 in
\cite{DuLG}. See also
Aldous \cite{Al2} and Duquesne \cite{Duq} for related statements.

Before stating our main results, we still need to introduce important
random measures
associated with stable trees. For every $a>0$, we can define
$\Theta_\alpha(d\t)$
a.e. a random measure
$\ell^a$ on the level set $\t(a):=\{\sigma\in\t:d(\rho,\sigma)=a\}$,
which is in a sense uniformly spread over that level set: For every
$\varepsilon>0$, write
$\t_\varepsilon(a)$ for the finite subset of $\t(a)$ consisting of
those vertices
which have descendants at level $a+\varepsilon$, then for every
bounded continuous
function $\varphi$ on $\t$,
$$\langle \ell^a,\varphi\rangle=\lim_{\varepsilon\da 0}{1\over v(\varepsilon)}
\sum_{\sigma\in \t_\varepsilon(a)}\varphi(\sigma).$$
We refer to Section 4.2
of \cite{DuLG} for the construction and main properties of these
``local time'' measures. The uniform measure $\m=\m_{(\t)}$ on the tree $\t$
is then defined by
\begin{equation}
\label{uniform-measure}
\m=\int_0^\infty da\,\ell^a\,.
\end{equation}

We start with the case $\alpha=2$ where we can identify the exact Hausdorff
measure function for the tree $\t$ and its level sets. The notation
$h-m$ stands for the
Hausdorff measure associated with the function $h$.

\begin{theorem}
\label{Brownian-Haus}
For every $r\in(0,1/2)$, set
$$h(r)=r^2\log\log{1\over r}.$$
There exists a positive constant $C_0$ such that
$\Theta_2$ a.e., for every Borel subset $A$ of $\t$,
$$h-m(A)=C_0\,\m(A).$$
\end{theorem}

According to this theorem, the measure $\m$ coincides with
a certain Hausdorff measure on $\t$. This justifies the fact
that $\m$ is called the {\it uniform} measure on the tree.

The law $\Theta^{(1)}_2$ of the CRT is informally defined by
$\Theta^{(1)}_2=\Theta_2(d\t\mid \m(\t)=1)$. More precisely,
the CRT is coded by a Brownian excursion conditioned to have duration $1$
(in the sense explained below in Section 3), whereas $\Theta_2$ is the law
of the tree coded by a Brownian excursion under the It\^o measure. Since the
excursion normalized to have duration $1$ and the It\^o measure are related
by simple scaling transformations, the following
corollary is an immediate consequence of Theorem \ref{Brownian-Haus}.

\begin{corollary}
\label{CRT}
Theorem \ref{Brownian-Haus} remains valid, with the same constant
$C_0$, if $\Theta_2$ is
replaced by the law $\Theta^{(1)}_2$ of the CRT.
\end{corollary}

Let us now discuss level sets. The next theorem shows that the local
time measure
$\ell^a$ coincides with a certain Hausdorff measure on the level set $\t(a)$.

\begin{theorem}
\label{Brownian-local-Haus}
For every $r\in(0,1/2)$, set
$$\wt h(r)=r\log\log{1\over r}.$$
There exists a positive constant $\wt C_0$ such that
for every $a>0$, one has
$\Theta_2$ a.e. for every Borel subset $A$ of $\t(a)$,
$$\wt h-m(A)=\wt C_0\,\ell^a(A).$$
\end{theorem}

When $1<\alpha<2$, we are unable to identify an exact Hausdorff
measure function
for the tree, but we still get rather precise information.

\begin{theorem}
\label{stable-Haus}
Suppose that $1<\alpha<2$.
For every $u\in\R$ and $r\in(0,e^{-1})$, set
$$h_u(r)=r^{\alpha\over \alpha-1}\,(\log{1\over r})^{1\over \alpha-1}
\,(\log\log {1\over r})^u.$$
Then,
\begin{description}
\item{\rm(i)} $h_u-m(\t)=\infty$ if $u>{1\over \alpha-1}$, $\Theta_\alpha$ a.e.
\item{\rm(ii)} $h_u-m(\t)=0$ if $u<0$, $\Theta_\alpha$ a.e.
\end{description}
\end{theorem}

The preceding results were announced, in a less precise form, in
Theorem 5.9 of \cite{DuLG}. Finally, we also have an analogue of Theorem
\ref{Brownian-local-Haus} in the stable case.

\begin{theorem}
\label{stable-local-Haus}
Suppose that $1<\alpha<2$ and let $a>0$.
For every $u\in\R$ and $r\in(0,e^{-1})$, set
$$\wt h_u(r)=r^{1\over \alpha-1}\,(\log{1\over r})^{1\over \alpha-1}
\,(\log\log {1\over r})^u.$$
Then,
\begin{description}
\item{\rm(i)} $\wt h_u-m(\t(a))=\infty$ if $u>{1\over \alpha-1}$,
$\Theta_\alpha$ a.e. on
$\{H(\t)>a\}$.
\item{\rm(ii)} $\wt h_u-m(\t(a))=0$ if $u<0$, $\Theta_\alpha$ a.e.
\end{description}
\end{theorem}

Let us briefly comment on the relation between these theorems and earlier
results. The Hausdorff dimension of stable trees was computed independently
in \cite{DuLG} and in \cite{HaM}. It is remarkable that the exact
Hausdorff measure
function of the tree under $\Theta_2$ (or of the CRT) is the same as
the one for a
transient Brownian path, which was
derived by Ciesielski and Taylor \cite{CT}
following earlier work of L\'evy. As we will see, some results from \cite{CT}
play a role in the proof of Theorem \ref{Brownian-Haus}. The preceding theorems
are also reminiscent of the very precise results about the Hausdorff
measure of the support and range of super-Brownian motion, which have been
obtained by Perkins and his co-authors (see \cite{Per1}, \cite{DIP},
\cite{LGP} and references therein). This should not come as a surprise since
superprocesses with a stable branching mechanism are easily constructed by
combining the genealogical structure of stable trees with independent
spatial motions (see e.g. Proposition 6.1 in \cite{DuLG}).

The paper is organized as follows. Section 2 gives the basic comparison results
for Hausdorff measures that are used in the proofs. Section 3
contains the proof
of Theorems \ref{Brownian-Haus} and \ref{Brownian-local-Haus}. Here we rely
on the coding of trees by Brownian excursions, which has been exploited in
other contexts, and in particular in the Brownian snake approach to
superprocesses \cite{Zurich}. Section 4 gives a few preliminary results about
stable trees, which are used in Section 5 to prove Theorems \ref{stable-Haus}
and \ref{stable-local-Haus}. In contrast with Section 3, we rely on general
properties of L\'evy trees that have been derived in \cite{DuLG}, and
in particular on the subtree decomposition along the ancestral line
of a typical
vertex (Theorem \ref{Palmdec} below). Section 5
also formulates conjectures for the exact Hausdorff measure of stable trees and
their level sets.

\section{Comparison results for Hausdorff measures}

In this section, we give a comparison result for Hausdorff
measures that will be used in the proofs below. For subsets of Euclidean space,
this result can be found as Lemmas 2 and 3 of Rogers and Taylor
\cite{RT} (see also
Theorem 1.4 in Perkins \cite{Per1} for a more precise formulation).
For the reader's convenience, and also because the arguments of \cite{RT} do
not extend immediately to the general setting which is considered here, we
provide a short proof below.

We consider a compact metric space $E$. For every $x\in E$ and $r>0$,
we denote by
$B(x,r)$ the open ball centered at $x$ with radius $r$. If
$c>1$ is fixed, we let ${\cal H}_c$ be the set of all monotone
increasing continuous functions
$g:[0,\infty)\to [0,\infty)$ such that $g(0)=0$ and $g(2r)\leq cg(r)$
for every $r\geq 0$.
As in the introduction, $g-m$ stands for the Hausdorff measure associated
with $g$. For any subset $A$ of $E$,
\be
\label{Hausdef}
g-m(A)=\lim_{\varepsilon\da 0}\Big(\inf_{(U_i)_{i\in I}\in{\cal
V}_\varepsilon(A)}
\sum_{i\in I} g({\rm diam}(U_i))\Big),
\ee
where ${\cal V}_\varepsilon(A)$ is the collection of all countable
coverings of $A$
by subsets of $E$
with diameter less than $\varepsilon$, and ${\rm diam}(U)$ denotes
the diameter of $U$.

\begin{lemma}
\label{comparison-Hausdorff}
Let $c>0$. There exist two positive constants $M_1$ and $M_2$ that depend
only on $c$, such that the following holds for every function $g\in{\cal H}_c$.
Let $\mu$ be a finite Borel measure on $E$ and let $A$ be a Borel subset of
$E$.
\begin{description}
\item{\rm (i)} If
$$\limsup_{n\to\infty} {\mu(B(x,2^{-n}))\over g(2^{-n})}\leq 1\
,\hbox{ for every }x\in A,$$
then,
$$g-m(A)\geq M_1\,\mu(A).$$
\item{\rm (ii)} If
$$\limsup_{n\to\infty} {\mu(B(x,2^{-n}))\over g(2^{-n})}\geq 1\
,\hbox{ for every }x\in A,$$
then,
$$g-m(A)\leq M_2\,\mu(A).$$
\end{description}
\end{lemma}

\proof (i) For every integer $n\geq 1$, set
$$A_n:=\{x\in A:\mu(B(x,2^{-k}))\leq 2\,g(2^{-k})\hbox{ for every }k\geq n\}.$$
By assumption, $A=\lim\ua A_n$ and so $\mu(A)=\lim\ua \mu(A_n)$. Now
fix $n\geq 1$
and consider a countable covering $(U_i)_{i\in I}$ of $A_n$ by sets
of positive diameter
strictly less than $2^{-n}$. For every $i\in I$, denote by $r_i>0$ the diameter
of $U_i$ and pick $x_i\in U_i\cap A_n$. Let $k_i\geq n$ be the unique
integer such
that $2^{-k_i-1}\leq r_i<2^{-k_i}$. Then, for every $i\in I$,
we have
$$U_i\cap A_n\subset \bar B(x_i,r_i)\subset B(x_i,2^{-k_i}).$$
Recalling the definition of $A_n$ it follows that
$$\sum_{i\in I} g(r_i)\geq c^{-1}\sum_{i\in I} g(2^{-k_i})\geq (2c)^{-1}
\sum_{i\in I} \mu(B(x_i,2^{-k_i}))\geq (2c)^{-1}\mu(A_n)$$
since the balls $B(x_i,2^{-k_i})$ cover $A_n$. From the definition of
Hausdorff measure we now get $g-m(A)\geq g-m(A_n)\geq (2c)^{-1}\mu(A_n)$
and the desired result follows by letting $n\ua \infty$.

\noindent(ii) From the general theory of Hausdorff measures (cf
Corollary 2, p.99
in \cite{Rog}), we know that
$$g-m(A)=\sup\{g-m(K):K\subset A,\,K\hbox{ compact}\}.$$
Hence we may assume in the proof that $A$ is compact.

Then let $\varepsilon>0$. By
assumption, for every $x\in A$, we may find $r_x\in(0,\varepsilon/8)$ such that
$$\mu(B(x,r_x))\geq {1\over 2}\,g(r_x).$$
By compactness, we may then find $x_1,\ldots,x_n\in A$, such that
$$A\subset \bigcup_{i=1}^n B(x_i,r_{x_i})$$
and we may assume that $r_{x_1}\geq r_{x_2}\geq \cdots\geq r_{x_n}$.
We can then
construct a finite subset $1=m_1<m_2<\cdots<m_\ell$ of
$\{1,2,\ldots,n\}$ in such a way
that if $y_j=x_{m_j}$ we have
$$A\subset \bigcup_{j=1}^\ell B(y_j,4r_{y_j})$$
and the balls $B(y_j,r_{y_j})$ and $B(y_{j'},r_{y_{j'}})$ are disjoint
if $j\not =j'$. In fact we start with $m_1=1$, and we proceed by
induction. Suppose that
we have constructed $m_1<m_2<\cdots<m_{p-1}$ in such a way that
$$\bigcup_{i=1}^{m_{p-1}} B(x_i,r_{x_i})\subset \bigcup_{j=1}^{p-1}
B(y_j,4r_{y_j}).$$
and the balls $B(y_j,r_{y_j})$, $1\leq j\leq p-1$ are disjoint. If
$$A\subset \bigcup_{i=1}^{p-1} B(y_j,4r_{y_j})$$
then the construction is complete. Otherwise we let $k>m_{p-1}$ be
the first integer
such that $B(x_k,r_{x_k})$ is not contained in the union of the
balls $B(y_j,4r_{y_j})$ for $j\leq p-1$, and we put $m_p=k$. Plainly,
$B(y_p,r_{y_p})\cap B(y_q,r_{y_q})=\varnothing$ if $1\leq q\leq p-1$, because
otherwise this would contradict the fact that the ball
$B(y_p,r_{y_p})=B(x_k,r_{x_k})$ contains
a point that does not belong to $B(y_q,4r_{y_q})$. This completes the
construction by induction.

\par Now the balls $B(y_j,4r_{y_j})$ provide a covering of $A$
by sets of diameter less than $\varepsilon$, and
$$\sum_{j=1}^\ell g(8r_{y_j})\leq c^3\sum_{j=1}^\ell g(r_{y_j})
\leq 2c^3 \sum_{j=1}^\ell \mu(B(y_j,r_{y_j}))\leq 2c^3\,\mu(A_\varepsilon)$$
where $A_\varepsilon$ stands for the $\varepsilon$-neighborhood of $A$.
Let $\varepsilon$ go to $0$ to get the desired result. \cq

\section{The Brownian tree}
In this section, we prove Theorem \ref{Brownian-Haus}
and Theorem \ref{Brownian-local-Haus}. We will make an extensive
use of the coding by Brownian excursions. Denote by ${\bf n}(de)$ the It\^o
measure of positive excursions of linear Brownian motion normalized so
that ${\bf n}(\sup e>\varepsilon)=\varepsilon^{-1}$, and by
$\zeta=\zeta(e)$ the duration of excursion
$e$. For every
$s,t\in[0,\zeta]$, we set
$$d_e(s,t)=e(s)+e(t)-2\,m_e(s,t)$$
where
$$m_e(s,t)=\inf_{s\wedge t\leq r\leq s\vee t} e(r).$$
We define an equivalence relation on $[0,\zeta]$ by
setting $s\sim t$ if $d_e(s,t)=0$. Then the quotient set
$\t_e:=[0,\zeta]/\sim$ equipped
with the metric $d_e$ is a random real tree (\cite{DuLG}
Theorem 2.1), whose root is by convention the equivalence class of $0$,
and the distribution of
$\t_e$ under
${\bf n}(de)$ is
$\Theta_2$. Furthermore, up to an unimportant multiplicative factor $2$
which we will ignore, the uniform measure
$\m$ on
$\t_e$ is just the image of  Lebesgue measure on $[0,\zeta]$ under
the canonical projection
from $[0,\zeta]$ onto $[0,\zeta]/\sim$, and similarly the local time measure
$\ell^a$ is the image of the usual Brownian local time measure at level
$a$. Therefore in
proving Theorem
\ref{Brownian-Haus} and Theorem \ref{Brownian-local-Haus}, we may and
will deal with the tree
$\t_e$ under ${\bf n}(de)$.

{\bf Proof of Theorem \ref{Brownian-Haus}.} We first establish the existence of
two positive constants $c_1$ and $c_2$ such that, ${\bf n}(de)$ a.e.
for every Borel subset $A$ of $\t_e$,
\be
\label{encad}
c_1\, \m(A)\leq h-m(A)\leq c_2\, \m(A).
\ee
{\bf Lower bound.} By abuse of notation we will often identify an element
$s$ of $[0,\zeta]$ with its equivalence class in
$\t_e=[0,\zeta]/\sim$. We first prove
that, ${\bf n}(de)$ a.e., for $\m$-almost all $s\in \t_e$,
one has
\be
\label{upperoccup}
\limsup_{\varepsilon\to 0} {\m(\{t\in \t_e:d_e(s,t)\leq
\varepsilon\})\over h(\varepsilon)}\leq C_1
\ee
for some finite constant $C_1$.

To prove (\ref{upperoccup}), we need a simple decomposition lemma
for the Brownian excursion. Assume that, on a certain probability space,
we are given two processes $(B_t,t\geq 0)$ and $(B'_t,t\geq 0)$
and for every $a\geq 0$ a probability measure $\Pi_a$ such that
$B$ and $B'$ are under $\Pi_a$ two independent Brownian motions
started at $a$. Also set
$$T=\inf\{t\geq 0: B_t=0\}\ ,\ T'=\inf\{t\geq 0: B'_t=0\}$$
and write $C(\R_+,\R)$ for the space of all continuous functions from
$\R_+$ into $\R$.

\begin{lemma}
\label{decomp}
For every nonnegative measurable function $F$
on $C(\R_+,\R)^2$,
$$\int {\bf n}(de) \int_0^\zeta ds F((e(s+t))_{t\geq 0},(e((s-t)_+))_{t\geq 0})
=2 \int_0^\infty da\,\Pi_a[F((B_{t\wedge T})_{t\geq 0},(B'_{t\wedge T'})_
{t\geq 0})].
$$
\end{lemma}

This is basically Bismut's decomposition of the Brownian excursion.
See \cite{LG1}, Lemma 1 for a simple proof (notice that our normalization of
It\^o's measure differs by a factor $2$ from the one in \cite{LG1}).

Note that by the definition of the distance $d_e$, and the preceding
identification of $\bf m$,
$$\m(\{t\in \t_e:d_e(s,t)\leq
\varepsilon\})
=\int_0^s dt\,1_{\{e(s)+e(t)-2m_e(s,t)\leq \varepsilon\}}
+\int_s^\zeta dt\,1_{\{e(s)+e(t)-2m_e(s,t)\leq \varepsilon\}}.$$
 From Lemma \ref{decomp}, we see that our claim (\ref{upperoccup})
will follow if we can prove that for every $a > 0$,
\be
\label{upperbis}
\limsup_{\varepsilon\to 0}
{\int_0^\infty dt\,1_{\{a+B_t-2I_t\leq \varepsilon\}}
+\int_0^\infty dt\,1_{\{a+B'_t-2I'_t\leq \varepsilon\}}\over
h(\varepsilon)}\leq C_1,\qquad\Pi_a\ {\rm a.s.}
\ee
where
$$I_t=\inf_{s\leq t} B_s\ ,\ I'_t=\inf_{s\leq t}B'_s.$$
By translation invariance, it is enough to consider the case $a=0$
in (\ref{upperbis}). A famous theorem of Pitman
states that the process $R_t:=B_t-2I_t$ is under $\Pi_0$ a
three-dimensional Bessel process started at $0$, that is, it has the
same distribution as the modulus of a three-dimensional
Brownian motion started from the origin. From estimates due
to Ciesielski and Taylor \cite{CT}, there exists a finite constant
$C_2$ such that
$$\limsup_{\varepsilon\to 0}{\int_0^\infty dt\,1_{\{R_t\leq \varepsilon\}}
\over h(\varepsilon)}=C_2,\qquad \Pi_0\ {\rm a.s.}$$
 From this and the analogous statement for $R'_t:=B'_t-2I'_t$,
we deduce (\ref{upperbis}), which completes the proof of (\ref{upperoccup}).
The lower bound in (\ref{encad}) then follows
from Lemma \ref{comparison-Hausdorff} (i).

\medskip
\noindent{\bf Upper bound.} From Lemma \ref{comparison-Hausdorff}
(ii), the upper
bound in (\ref{encad}) will follow if we can prove the existence of a
constant $K_1>0$ such
that,
${\bf n}(de)$ a.e.,
\be
\label{loweroccup}
h-m(\{s\in\t_e: \limsup_{\varepsilon\to 0} {\m(\{t\in
\t_e:d_e(s,t)\leq
\varepsilon\})\over h(\varepsilon)}\leq K_1\})=0.
\ee
For every integer $n\geq 0$, set $\varepsilon_n=2^{-n}$. We will
prove the existence of a constant $K_2$
such that, for every integer $n_0\geq 0$, ${\bf n}(de)$ a.e.,
\be
\label{loweroccup2}
h-m(\{s\in\t_e: e(s)>2^{-n_0}\hbox{ and }\m(\{t\in
\t_e:d_e(s,t)\leq
\varepsilon_p\})\leq K_2 h(\varepsilon_p),
\forall p\geq n_0\})=0.
\ee
Clearly, (\ref{loweroccup}) follows from (\ref{loweroccup2}).
To prove (\ref{loweroccup2}), we will need to introduce suitable
coverings of the sets
$$F^A_{n_0}:=\{s\in\t_e: 2^{-n_0}<e(s)<a\hbox{ and }\m(\{t\in
\t_e:d_e(s,t)\leq
\varepsilon_p\})\leq K_2 h(\varepsilon_p),
\forall p\geq n_0\},$$
where $A$ is a positive integer. For every $n\geq 0$, consider the
sequence of stopping times defined inductively as follows
$$T^n_0=0,\ T^n_1=\inf\{s\geq 0:e(s)=2^{-n}\},\
T^n_{k+1}=\inf\{s\geq T^n_k:|e(s)-e(T^n_k)|=2^{-n}\},$$
where $\inf\varnothing=\infty$ by convention. The sequence
$(2^ne(T^n_k)1_{\{T^n_k<\infty\}},k\geq 0)$ is distributed under
${\bf n}(de\mid T^n_1<\infty)$ as a positive excursion of simple
random walk. In
particular, for every integer $j\geq 1$,
\be
\label{localrw}
{\bf n}\Big(\sum_{k=0}^\infty
1_{\{T^n_k<\infty,\,e(T^n_k)=j2^{-n}\}}\Big)= 2\,{\bf n}(T^n_1<\infty)=
2^{n+1}.
\ee
Let $s\in F^A_{n_0}$ and $n\geq n_0+2$. There exists a unique
integer $k>0$ such that $s\in[T^n_k,T^n_{k+1})$. From our definitions,
we have then
$$d_e(T^n_k,s)\leq 3\cdot 2^{-n}.$$
As a consequence, for every $p\in\{n_0,n_0+1,\ldots,n-2\}$, we have
$$\{t\in
\t_e:d_e(s,t)\leq
\varepsilon_p\}\supset\{t\in
\t_e:d_e(T^n_k,t)\leq
\varepsilon_p/4\}.$$
It follows that
\be
\label{cover1}
F^A_{n_0}
\subset \bigcup_{k\in I_{n_0,n}}[T^n_k,T^n_{k+1})
\ee
where
\ba &&I_{n_0,n}=
\{k\geq 0:T^n_k<\infty,\ 2^{-n_0}\leq e(T^n_k)\leq A\\
&&\hspace{4cm}\hbox{ and }\int_{T^n_k}^\zeta dt\,1_{\{d_e(T^n_k,t)\leq
\varepsilon_p/4\}}
\leq K_2h(\varepsilon_p),
\ \forall p\in\{n_0,\ldots,n-2\}\}.
\ea
To bound the cardinality $\# I_{n_0,n}$ of the set $I_{n_0,n}$,
we use the strong Markov property under the excursion measure
to write, for every $k\geq 1$,
\ba &&{\bf n}(k\in I_{n_0,n})
={\bf n}\Big({\bf 1}_{\{T^n_k<\infty,\,2^{-n_0}\leq e(T^n_k)\leq A\}}\\
&&\hspace{3cm}\times \,\Pi_{e(T^n_k)}\Big[\int_0^T
dt\,1_{\{e(T^n_k)+B_t-2I_t\leq
\varepsilon_p/4\}}\leq K_2h(\varepsilon_p),
\ \forall p\in\{n_0,\ldots,n-2\}\Big]\Big).
\ea
Using again Pitman's theorem recalled above, we have
for every $a\geq 2^{-n_0}$,
\ba&&\Pi_{a}\Big[\int_0^T dt\,1_{\{a+B_t-2I_t\leq
\varepsilon_p/4\}}\leq K_2h(\varepsilon_p),
\ \forall p\in\{n_0,\ldots,n-2\}\Big]\\
&&\quad=\Pi_0\Big[\int_0^\infty dt\,1_{\{R_t\leq \varepsilon_p/4\}}
\leq K_2h(\varepsilon_p),\ \forall p\in\{n_0,\ldots,n-2\}\Big]
\ea
where $R$ is under $\Pi_0$ a three-dimensional Bessel process started
at $0$. It follows from Theorem 1.2 in \cite{LG0} that, provided
$n_0$ is large enough, we can choose $K_2$ sufficiently small
so that the last probability is bounded above by
$$\exp(-c(n-n_0)^{1/2})$$
for some positive constant $c$. Hence,
$${\bf n}(k\in I_{n_0,n})\leq \exp(-c(n-n_0)^{1/2})\,
{\bf n}(T^n_k<\infty,\ 2^{-n_0}\leq e(T^n_k)\leq A)$$
and by summing over $k$, and using (\ref{localrw})
$${\bf n}(|I_{n_0,n}|)\leq \exp(-c(n-n_0)^{1/2})\,A\,2^{2n+1}.$$
In particular, by Fatou's lemma,
$$\liminf_{n\to \infty}
2^{-2n}\,\exp(c(n-n_0)^{1/2})\,|I_{n_0,n}|<\infty$$
${\bf n}(de)$ a.e. Now recall (\ref{cover1}) and note that the diameter
(with respect to the distance $d_e$) of each interval
$[T^n_k,T^n_{k+1}]$ is bounded above by $4\,2^{-n}$. Our claim
(\ref{loweroccup2}) then follows from the definition of Hausdorff
measures. This completes the proof of (\ref{encad}).

\medskip
Theorem \ref{Brownian-Haus} can be deduced from the bounds (\ref{encad})
and an appropriate zero-one law. This is similar to the argument used
in Section 7 of \cite{LGP}, but there are some differences.

Let us write $p_e$ for the canonical projection from
$[0,\zeta]$ onto $\t_e=[0,\zeta]/\sim$.  We first observe
that, for every $0\leq s\leq t\leq \zeta$, the quantity
$h-m(p_e([s,t]))$ is a mesurable function of $e$. To see this, note that
in the definition (\ref{Hausdef}) of $h-m(p_e([s,t]))$, we may
restrict our attention to finite coverings with balls (use compactness and
the fact that any subset
of a real tree is contained in a closed ball with the same diameter). Moreover,
it is enough to consider balls with rational diameter, and with a center of
the form $p_e(r)$ for some rational number $r\in[0,\zeta]$. The
desired measurability
property then follows easily.

  We then define a finite measure
$\nu$ on $[0,\zeta]$ by setting, for every $t\in[0,\zeta]$,
$$\nu([0,t]) = h-m( p_e([0,t]) ).$$
Plainly, the mapping $t\la h-m(p_e([0,t]))$ is continuous and so $\nu$
is nonatomic.
Then we have also, for every $0\leq s\leq t\leq \zeta$,
$$\nu([s,t]) = h-m( p_e([s,t]) ).$$
Indeed, this is a consequence of the following observation: If $0\leq
u<v\leq s<t\leq\zeta$,
the set $p_e([u,v])\cap p_e([s,t])$ is contained in the ancestral line
of $p_e(s)$, and so we must have
$$h-m(p_e([u,v])\cap p_e([s,t]))=0.$$
Since $\m$ is obtained as
the image of Lebesgue measure under $p_e$, it is easy to verify that
$\m(p_e([s,t]))=t-s$ for every $0\leq s\leq t\leq \zeta$, ${\bf n}(de)$ a.e.
 From the bounds (\ref{encad}), we get ${\bf n}(de)$ a.e.
for every $0\leq s\leq t\leq \zeta$,
$$ c_1 (t-s) \leq \nu([s,t])  \leq c_2 (t-s).$$
Hence the measure $\nu$ is absolutely continuous with respect to
Lebesgue measure on $[0,\zeta]$, and by a standard differentiation
theorem its density is equal almost everywhere to
\be
\label{zero-one00}
{d\nu\over dt}=\lim_{\epsilon \to 0}  {\nu([t-\epsilon,t+\epsilon])\over
2\epsilon} = \lim_{\epsilon \to 0} {h-m( p_e
([t-\epsilon,t+\epsilon]))\over 2\epsilon}.
\ee
It is easy to see that the quantity $h-m( p_e ([t-\epsilon,t+\epsilon]))$
is a measurable function of the path $(e(t+u)-e(t),-\epsilon\leq
u\leq \epsilon)$.
Hence we can use Lemma \ref{decomp} and the standard $0-1$ law for Brownian
motion to get that the last limit
in (\ref{zero-one00}) must be equal to a constant $C_0\in[0,\infty]$,
$dt$ a.e., ${\bf n}(de)$ a.e. Obviously, $C_0\in[c_1,c_2]$ and in
particular $0<C_0<\infty$.

We have thus $h-m(A) = C_0\, \m(A)$ for every subset of the tree
of the form $p_e([s,t])$, or for any finite union of such sets.
However, every open subset $U$ of the tree is the increasing limit
of a sequence of such unions (note that $p_e^{-1}(U)$
is a countable union of open intervals). Hence $h-m(A) = C_0\, \m(A)$
for every open subset $A$ of $\t_e$, which is enough to complete the proof.
  \cq

\bigskip

{\bf Proof of Theorem \ref{Brownian-local-Haus}.} We now turn to the
Hausdorff measure of level
sets of
$\t_e$. Recall that
$$\t_e(a)=\{s\in\t_e:d_e(0,s)=a\}=\{s\in \t_e:e(s)=a\}.$$
We also denote by $(\ell^a_s,s\geq 0)$ the (Brownian) local time process
of $e$ at level $a$. Then the measure $\ell^a(ds)$ associated with the
increasing function
$s\to\ell^a_s$ can be interpreted as a
measure on $\t_e(a)$ and indeed coincides with the one discussed in the
introduction (up to a multiplicative factor $2$ which is irrelevant
for our purposes).
Moreover, for every nonnegative measurable function $F$
on $C(\R_+,\R)^2$,
\be
\label{Bismut-local}
\int {\bf n}(de) \int_0^\zeta \ell^a(ds)\, F((e(s+t))_{t\geq
0},(e((s-t)_+))_{t\geq 0})
=2\,\Pi_a[F((B_{t\wedge T})_{t\geq 0},(B'_{t\wedge T'})_
{t\geq 0})].
\ee
This formula is easily derived from Lemma \ref{decomp} and the usual
approximations
of Brownian local time.

As in the proof
of Theorem \ref{Brownian-Haus}, we first establish the existence of
two positive constants
$\wt c_1$ and $\wt c_2$ such that, ${\bf n}(de)$ a.e. for every Borel
subset $A$ of $\t(a)$,
\be
\label{encadbis}
\wt c_1\, \ell^a(A)\leq \wt h-m(A)\leq \wt c_2\,\ell^a(A).
\ee
{\bf Lower bound.} Similarly as in the proof
of the lower bound in (\ref{encad}), it is enough to show that there
exists a constant $C'_1$ such that,
${\bf n}(de)$ a.e., for $\ell^a$-almost all $s\in \t_e(a)$,
\be
\label{upperoccuplevel}
\limsup_{\varepsilon\to 0} {\ell^a(\{t\in \t_e(a):d_e(s,t)\leq
\varepsilon\})\over \wt h(\varepsilon)}\leq C'_1.
\ee
If $0<\varepsilon<a$ and $s,t\in\t_e(a)$, we have $d_e(s,t)\leq \varepsilon$
if and only if $m_e(s,t)\leq \varepsilon/2$. From this observation
and (\ref{Bismut-local}), we see that (\ref{upperoccuplevel})
will follow if we can verify that, for every $a>0$, $\Pi_a$ a.s.,
\be
\label{leveltech}
\limsup_{\varepsilon\to
0}{L^a_{T_{a-\varepsilon}}(B)+L^a_{T'_{a-\varepsilon}}(B')
\over \wt h(\varepsilon)}\leq C'_1,
\ee
where $(L^a_t(B),t\geq 0)$ is the local time process of $B$ at level $a$,
and $T_{a-\varepsilon}=\inf\{t:B_t=a-\varepsilon\}$, with a similar
notation for $L^a_t(B')$ and $T'_{a-\varepsilon}$.

It is well known that the distribution of $L^a_{T_{a-\varepsilon}}(B)$
under $\Pi_a$ is exponential with mean $2\varepsilon$. Therefore an
application of the Borel-Cantelli lemma immediately shows that,
for $\varepsilon_n=2^{-n}$,
$$\limsup_{n\to\infty}{L^a_{T_{a-\varepsilon_n}}(B)\over \wt h(\varepsilon_n)}
\leq 1.$$
It readily follows that (\ref{leveltech}) holds with $C'_1=8$.

\bigskip
\noindent{\bf Upper bound.} This is similar to the proof of the upper bound
in (\ref{encad}). It is enough to show that there is
a constant $K'$, not depending on $a$, such that, ${\bf n}(de)$ a.e.,
\be
\label{loweroccupbis}
\wt h-m(\{s\in\t_e(a): \limsup_{\varepsilon\to 0} {\ell^a(\{t\in
\t_e(a):d_e(s,t)\leq
\varepsilon\})\over \wt h(\varepsilon)}\leq K' \})=0.
\ee
This requires finding good coverings for the sets
$$G_{n_0}=\{s\in \t_e(a):\;\ell^a(\{t\in\t_e(a):d_e(s,t)\leq \varepsilon_p\})
\leq K'\wt h(\varepsilon_p),\,\forall p\geq n_0\},$$
for all $n_0$ sufficiently large. Fix $a>0$ and $n_0\geq 1$
such that $2^{-n_0}<a$. To cover $G_{n_0}$, introduce the stopping
times defined for every $n\geq n_0$,
$$T^n_0=\inf\{s:e(s)=a\}\ ,\ T^n_1=\inf\{s>T^n_0:|e(s)-a|=2^{-n}\}$$
and by induction,
$$T^n_{2k}=\inf\{s>T^n_{2k-1}:e(s)=a\}\ ,\
T^n_{2k+1}=\inf\{s>T^n_{2k}:|e(s)-a|=2^{-n}\}.$$
It is easy to verify that
\be
\label{leveltech2}
{\bf n}\Big(\sum_{k=0}^\infty 1_{\{T^n_{2k}<\infty\}}\Big)\leq C_{(a)}2^n
\ee
where the constant $C_{(a)}$ only depends on $a$.

In a way very similar to the proof of Theorem \ref{Brownian-Haus}, we have
\be
\label{leveltech3}
G_{n_0}\subset \bigcup_{k\in
J_{n_0,n}}[T^n_{2k},T^n_{2k+1}]
\ee
where
$$J_{n_0,n}=\{k:T^n_{2k}<\infty,\;\int_{T^n_{2k}}^\zeta
d\ell^a_t\,1_{\{d_e(T^n_{2k},t)\leq \varepsilon_p/4\}}\leq
K'\wt h(\varepsilon_p),\,\forall p\in\{n_0,\ldots,n-2\}\}.$$
By the strong Markov property at time $T^n_{2k}$,
$${\bf n}(k\in J_{n_0,n})
={\bf n}\Big(T^n_{2k}<\infty,\,\Pi_a[L^a_{T_{a-\varepsilon_p/8}}\leq
K'\wt h(\varepsilon_p),\forall p\in\{n_0,\ldots,n-2\}]\Big).$$
Now note that the variables $L^a_{T_{a-\varepsilon_p/8}}
-L^a_{T_{a-\varepsilon_p/16}}$, $p\in\{n_0,\ldots,n-2\}$
are independent under $\Pi_a$. Moreover, conditionally on the event that
it is strictly positive, which has probability $1/2$, the variable
$L^a_{T_{a-\varepsilon_p/8}}
-L^a_{T_{a-\varepsilon_p/16}}$
is exponentially distributed with mean $\varepsilon_p/4$.
It follows that
\ba
&&\Pi_a[L^a_{T_{a-\varepsilon_p/8}}\leq
K'\wt h(\varepsilon_p),\forall p\in\{n_0,\ldots,n-2\}]\\
&&\quad\leq \Pi_a[L^a_{T_{a-\varepsilon_p/8}}-L^a_{T_{a-\varepsilon_p/16}}\leq
K'\wt h(\varepsilon_p),\forall p\in\{n_0,\ldots,n-2\}]\\
&&\quad=\prod_{p=n_0}^{n-2}(1-{1\over 2}\exp(-4K'\log\log 2^p)).
\ea
If $K'<1/8$, the latter quantity is bounded above by
$\exp(-n^{1/2})$ for $n$ large. Therefore we get for all
$n$ sufficiently large,
$${\bf n}(k\in J_{n_0,n})\leq \exp(-n^{1/2})\,{\bf n}(T^n_{2k}<\infty).$$
By combining this with (\ref{leveltech2}), and using Fatou's lemma,
we arrive at
\be
\label{endpr}
\liminf_{n\to\infty} 2^{-n}\exp(n^{1/2})\,\#J_{n_0,n}<\infty,
\ee
${\bf n}(de)$ a.e.
Since by construction the $d_e$-diameter of each interval
$[T^n_{2k},T^n_{2k+1}]$ is bounded above by $4\,2^{-n}$, (\ref{endpr}) and
(\ref{leveltech3}) lead to $\wt h-m(G_{n_0})=0$, which
completes the proof of (\ref{loweroccupbis}) and of
the bounds (\ref{encadbis}).

\medskip
The end of the proof is now similar to the final part of
the proof of Theorem \ref{Brownian-Haus}. We introduce the random
measure $\wt\nu$ on $[0,\zeta]$ defined by
$$\wt\nu([0,t])=\wt h-m(p_e([0,t])\cap \t(a)).$$
The bounds (\ref{encadbis}) imply that $\wt\nu(dt)$ is absolutely
continuous with respect to $\ell^a(dt)$, and that its density
is bounded below and above by $\wt c_1$ and $\wt c_2$ respectively.
A zero-one law argument, now relying on (\ref{Bismut-local}), shows that
this density is equal to a constant $\wt C_0$, $\ell^a(dt)$ a.e.,
${\bf n}(de)$ a.e. Moreover this constant does not depend on $a$. We leave
details to the reader.
  \cq

\section{Preliminaries about stable trees}

In this section we collect the basic facts about stable trees
that will be needed in the proof of Theorem \ref{stable-Haus}
and Theorem \ref{stable-local-Haus}.
We refer to \cite{DuLG} for additional details.

We fix $\alpha\in(1,2)$.
As in the introduction above, we write $\Theta_\alpha$ for the
distribution of the stable tree with index $\alpha$. In the terminology
of \cite{DuLG}, this corresponds to the measure $\Theta$ associated
with the branching mechanism function $\psi(u)=u^\alpha$.
Note that $\Theta_\alpha$ is a $\sigma$-finite measure on the
space $\T$ of (rooted) $\R$-trees, which puts no mass on the trivial
tree consisting only of the root.

In the same way as in the previous section, $\Theta_\alpha$ can be
defined and studied in terms of its coding function. However,
the role of the Brownian excursion in the case $\alpha=2$ is now
played by the stable height process, which is a less tractable
probabilistic object. For this reason, rather than using the coding
function as we
did in the case $\alpha=2$, we will state here the key properties
of the stable tree that are relevant to our study, and that
can be found in \cite{DuLG}.

We already mentioned the scaling invariance property of $\Theta_\alpha$:
For every $r>0$, the distribution of the scaled tree $r\t$
under $\Theta_\alpha$ is $r^{{1\over \alpha-1}}\Theta_\alpha$. We can
also express the local times
$\ell^a_{(r\t)}$ and uniform measure $\m_{(r\t)}$ of the scaled tree
$r\t$ in terms of the local times $\ell^a_{(\t)}$ and
uniform measure $\m_{(\t)}$ of the
tree $\t$. Precisely, considering only the total masses of these
random measures, we have $\Theta_\alpha$ a.e.,
\be
\label{scale-local}
\langle \ell^a_{(r\t)},1\rangle =r^{{1\over
\alpha-1}}\langle\ell^{a/r}_{(\t)},1
\rangle\ ,\ \langle \m_{(r\t)},1\rangle=r^{{\alpha\over \alpha-1}}\langle
\m_{(\t)},1\rangle.
\ee
This can be checked from the approximation of
local time recalled in the introduction above.

Informally, the tree $\t$ under $\Theta_\alpha$ describes the
genealogy of descendants of a single individual in a continuous-state
branching process with branching mechanism $\psi(u)=u^\alpha$.
The total mass $\langle \ell^a,1\rangle$ then corresponds to the
population at time (or level) $a$. To make this more precise, we can state
the following ``Ray-Knight property'' of local times. Let $x>0$ and let
$$\sum_{i\in I} \delta_{\t_i}$$
be  a Poisson point
measure on $\T$ with intensity $x\Theta_\alpha$. The real-valued
process $(X_t)_{t\geq 0}$ defined by
$$\left\{
\begin{array}{ll}
X_t={\displaystyle\sum_{i\in I} \langle \ell^t_{(\t_i)},1\rangle}\quad
&\hbox{if } t>0\;,\\
X_0=x&\\
\end{array}
\right.
$$
is a continuous-state
branching process with branching mechanism $\psi(u)=u^\alpha$,
started at $X_0=x$. This means that $(X_t)_{t\geq 0}$ is a Feller
Markov process on $\R_+$ and that the Laplace transform of its
semigroup is determined as follows: For every $\lambda>0$,
$$E[\exp(-\lambda X_t)\mid X_0=x]=\exp(-x\,u_t(\lambda))$$
where $(u_t(\lambda))_{t\geq 0}$ is determined from the integral
equation
$$u_t(\lambda)+\int_0^t ds \,u_s(\lambda)^\alpha = \lambda,$$
so that
\be
\label{Laplexp}
u_t(\lambda)=(\lambda^{1-\alpha}+(\alpha-1)t)^{1\over 1-\alpha}.
\ee
Note that we have also
$$u_t(\lambda)=\Theta_\alpha(1-\exp-\lambda\langle \ell^t,1\rangle)$$
from the exponential formula for Poisson measures. Using the Markov property
of $X$, one easily derives similar integral equations for
finite-dimensional marginal distributions of $(X_t)_{t\geq 0}$: See e.g.
Section II.3  of \cite{Zurich} where the more general setting of
superprocesses is considered. We will need the following particular
case: For every $\gamma,\lambda>0$,
the function
$$v_t(\gamma,\lambda)=\Theta_\alpha\Big(1-\exp\Big(-\gamma\int_0^tds\,\langle
\ell^s,1\rangle-\lambda\langle \ell^t,1\rangle\Big)\Big)$$
solves the integral equation
$$v_t+\int_0^t ds\,(v_s)^\alpha =\gamma t + \lambda.$$
Recall that $\t(r)$ denotes the level set of $\t$ at level $r$, and $H(\t)$
stands for the height of $\t$. We
also use the notation $\t_{\leq r}$ for the set
$\{\sigma\in\t:d(\rho,\sigma)\leq r\}$.

\begin{lemma}
\label{technical-stable}
There exist two positive constants $c_\alpha$ and $C_\alpha$ such that,
for every $b>0$,
$$\Theta_\alpha\Big(\m(\t_{\leq 1})\leq b\;,\;H(\t)\geq 1\Big)
\leq C_\alpha\,\exp(-c_\alpha b^{-{\alpha-1\over \alpha}}).$$
\end{lemma}

\proof
With the preceding notation, set
$$v^0_t(\gamma)=v_t(\gamma,0)=
\Theta_\alpha\Big(1-\exp\Big(-\gamma\int_0^tds\,\langle
\ell^s,1\rangle\Big)\Big)$$
so that $v^0_t=v^0_t(\gamma)$ solves the integral equation
$$v^0_t+\int_0^t ds\,(v^0_s)^\alpha =\gamma t .$$
It follows, that, for every $t\geq 0$, $v^0_t(\gamma)\in[0,\gamma^{1/\alpha})$
is determined by
\be
\label{eq-v-0}
\int_0^{v^0_t(\gamma)} {dy\over \gamma-y^\alpha}=t.
\ee
Similarly, if $v^\infty_t(\gamma)=\lim\ua v_t(\gamma,\lambda)$
as $\lambda\ua\infty$, we have
$$v^\infty_t(\gamma)=
\Theta_\alpha\Big(1-{\bf
1}_{\{\ell^t=0\}}\exp\Big(-\gamma\int_0^tds\,\langle
\ell^s,1\rangle\Big)\Big)$$
and $v^\infty_t(\gamma)\in(\gamma^{1/\alpha},\infty)$
is determined from the equation
$$\int_{v^\infty_t(\gamma)}^\infty {dy\over y^\alpha -\gamma}=t.$$
Simple analytic arguments show that
$$v^0_t(\gamma)=\gamma^{1/\alpha}(1-e^{-\varphi_t(\gamma)})$$
where
$$\lim_{\gamma\to\infty} \gamma^{{1\over
\alpha}-1}\varphi_t(\gamma)=\alpha t.$$
Similarly,
$$v^\infty_t(\gamma)=\gamma^{1/\alpha}(1+e^{-\phi_t(\gamma)})$$
where
$$\lim_{\gamma\to\infty} \gamma^{{1\over
\alpha}-1}\phi_t(\gamma)=\alpha t.$$
Now observe that
$$\Theta_\alpha\Big({\bf
1}_{\{\ell^1\not =0\}}\exp\Big(-\gamma\int_0^1ds\,\langle
\ell^s,1\rangle\Big)\Big)=v^\infty_1(\gamma)-v^0_1(\gamma)=
\gamma^{1/\alpha}(e^{-\varphi_1(\gamma)} + e^{-\phi_1(\gamma)}).$$
Furthermore, by construction, ${\bf m}(\t_{\leq 1})=\int_0^1
ds\,\langle
\ell^s,1\rangle$ and $\{\ell^1\not =0\}=\{H(\t)\geq 1\}$,
$\Theta_\alpha$ a.e.
(cf Theorem 4.2 in \cite{DuLG}). We get
$$\lim_{\gamma\to\infty}
{\log \Theta_\alpha({\bf 1}_{\{H(\t)\geq 1\}}\,\exp(-\gamma
{\bf m}(\t_{\leq 1})))\over \gamma^{1-{1\over \alpha}}}
=-\alpha.$$
The estimate of the lemma now follows. \cq

\smallskip

An important role in the next section will be played by
a subtree decomposition along the ancestral line of a
randomly chosen vertex. For the reader's convenience, we now recall this
result.

We first introduce the relevant notation. Let $\t\in \T$ and $\sigma\in\t$.
Denote by
$\llbracket \rho(\t),\sigma\rrbracket$ the line segment
from the root $\rho$ to $\sigma$, that is the ancestral line of $\sigma$.
If $\sigma,\sigma'\in\t$, the notation $\sigma\wedge \sigma'$ stands
for the most recent
common ancestor to $\sigma$ and $\sigma'$ (equivalently,
$\llbracket \rho,\sigma\rrbracket\cap \llbracket \rho,\sigma'\rrbracket
=\llbracket \rho,\sigma\wedge \sigma'\rrbracket$).
Denote by
$\t^{(j),\circ}$,
$j\in{\cal J}$ the connected components of
the open set $\t\backslash \llbracket \rho,\sigma\rrbracket$, and note
that for every
$j\in {\cal J}$, $\sigma_j:=\sigma\wedge \tau$ does not depend on the
choice of
$\tau\in \t^{(j),\circ}$.  Furthermore,
$\t^{(j)}:=\t^{(j),\circ}\cup\{\sigma_j\}$ is a (compact rooted) $\R$-tree
with root $\sigma_j$. The trees $\t^{(j)}$, $j\in {\cal J}$ can
be interpreted as the subtrees
of $\t$ originating from the segment
$\llbracket \rho,\sigma\rrbracket$. We put
$${\cal M}_\sigma=\sum_{j\in {\cal J}}
\delta_{(d(\rho(\t),\sigma_j),\t^{(j)})},$$ thus defining a point measure
on $[0,\infty)\times \T$.

\begin{theorem}
\label{Palmdec}
For every $a>0$ and
every nonnegative measurable function $\Phi$ on $[0,\infty)\times \T$,
$$\Theta_\alpha\Big(\int \ell^a(d\sigma)\,\exp-\langle
{\cal M}_\sigma,\Phi\rangle\Big) =\exp\Big(-\alpha\int_0^a
dt\,\Big(\Theta_\alpha(1-\exp-\Phi(t,\cdot))\Big)^{\alpha-1}\Big).$$
\end{theorem}

This is the case $\psi(u)=u^\alpha$ in Theorem 4.5 of \cite{DuLG}.

\section{The Hausdorff measure of the stable tree}

In this section we prove Theorem \ref{stable-Haus}
and Theorem \ref{stable-local-Haus}.
We keep the notation and assumptions of the preceding section.

{\bf Proof of Theorem \ref{stable-Haus}}. Part (i)
is an immediate consequence of the following proposition.

\begin{proposition}
\label{lower-Haus}
Suppose that $h:[0,\infty)\to[0,\infty)$ is a monotone increasing
function that can be written in the form $h(r)=r^{\alpha\over
\alpha-1}g(r)$ where $g$ is monotone decreasing in a neighborhood of
the origin. Then the condition
\be
\label{sum-condi}
\sum_{n=1}^\infty g(2^{-n})^{-(\alpha-1)}<\infty
\ee
implies that $h-m(\t)=\infty$, $\Theta_\alpha$ a.e.
\end{proposition}

\noindent{\bf Conjecture.} If (\ref{sum-condi}) fails, then
$h-m(\t)=0$, $\Theta_\alpha$
a.e.

\smallskip
\proof Let $a>0$. If $\sigma\in\t$ and $\varepsilon>0$, denote by
$B(\sigma,\varepsilon)$ the closed ball of radius $\varepsilon$ centered
at $\sigma$.
Then Theorem
\ref{Palmdec} implies that, for every
$\lambda>0$ and $\varepsilon\in(0,a]$,
\be
\label{lower-tech1}
\Theta_\alpha\Big(\int \ell^a(d\sigma)\,\exp-\lambda{\bf
m}(B(\sigma,\varepsilon))\Big)
=\exp\Big(-\alpha\int_{a-\varepsilon}^a \Phi_{\varepsilon,\lambda,a}(t)
^{\alpha-1}\,dt\Big),
\ee
where for $t\in[a-\varepsilon,a]$,
$$\Phi_{\varepsilon,\lambda,a}(t)=
\Theta_\alpha\Big(1-\exp-\lambda\m(B(\rho,\varepsilon-(a-t)))\Big).$$
Details of the derivation of (\ref{lower-tech1}) can be found on p.593-594
of \cite{DuLG}, where a similar formula is derived in greater
generality.

In agreement with the notation of the preceding section, we put
$$v^0_r(\lambda)=\Theta_\alpha(1-\exp-\lambda \m(B(\rho,r)))
=\Theta_\alpha(1-\exp-\lambda \m(\t_{\leq r})).$$
Then, (\ref{lower-tech1}) can be rewritten in the form
\be
\label{lower-tech3}
\Theta_\alpha\Big(\int \ell^a(d\sigma)\,\exp-\lambda{\bf
m}(B(\sigma,\varepsilon))\Big)
=\exp\Big(-\alpha\int_0^\varepsilon dr\,v^0_r(\lambda)^{\alpha-1}\Big).
\ee

Recall that $v^0_r(\lambda)$ is determined from equation
(\ref{eq-v-0}). From this
equation one immediately derives the following scaling property: For
every $\varepsilon>0$,
\be
\label{scaling-v}
v^0_{\varepsilon r}(\lambda)=\varepsilon^{-{1\over
\alpha-1}}\,v^0_r(\varepsilon^{\alpha\over \alpha-1}\lambda).
\ee
Hence,
$$\int_0^\varepsilon dr\,v^0_r(\lambda)^{\alpha-1}
=\int_0^1 dr\,v^0_r(\varepsilon^{\alpha\over
\alpha-1}\lambda)^{\alpha-1}.$$
If we substitute this identity into (\ref{lower-tech3}) and
replace $\lambda$ by $\varepsilon^{-{\alpha\over \alpha-1}}\lambda$,
we arrive at
\be
\label{lower-tech4}
\Theta_\alpha\Big(\int \ell^a(d\sigma)\,\exp(-
\lambda \varepsilon^{-{\alpha\over
\alpha-1}} {\bf m}(B(\sigma,\varepsilon)))\Big)
=\exp\Big(-\alpha \int_0^1 dr\,v^0_r(\lambda)^{\alpha-1}\Big).
\ee
The local time measure $\ell^a$ satisfies
$\Theta_\alpha(\langle\ell^a,1\rangle)=1$ (take $\lambda=0$ in
(\ref{lower-tech4})).
Thus, $\Theta_\alpha(d\t)\ell^a(d\sigma)$ defines a probability measure on the
set of ``pointed $\R$-trees'', that is pairs consisting of an
$\R$-tree $\t$  and a distinguished point $\sigma\in\t$ (in addition to
the root).
Denote by $\mu$ the law of $\varepsilon^{-{\alpha\over
\alpha-1}} {\bf m}(B(\sigma,\varepsilon))$ under the probability
measure $\Theta_\alpha(d\t)\ell^a(d\sigma)$. By (\ref{lower-tech4}),
this law does not depend on the choice of $\varepsilon$ and $a$,
provided that $0<\varepsilon\leq a$. Furthermore, the Laplace transform
of $\mu$ is given by
$$\int \mu(dx)\,e^{-\lambda x}=\exp\Big(-\alpha\int_0^1
dr\,v^0_r(\lambda)^{\alpha-1}\Big).$$
By monotone convergence,
$${v^0_r(\lambda)\over \lambda}=\Theta\Big({1-\exp(-\lambda \m(B(\rho,r)))
\over \lambda}\Big)\build{\la}_{\lambda\da 0}^{}
\Theta(\m(B(\rho,r)))=r$$
and
$$\lambda^{1-\alpha} \int_0^1 dr\,v^0_r(\lambda)^{\alpha-1}
\build{\la}_{\lambda\da 0}^{}
\int_0^1 dr\,r^{\alpha-1}={1\over \alpha}.$$
It follows that
$$\int \mu(dx)\,e^{-\lambda
x}=1-\lambda^{\alpha-1}+o(\lambda^{\alpha-1})$$ as $\lambda\to 0$.
Consequently, there exists a constant $C$ such that, for every $y>0$,
\be
\label{lower-tech5}
\mu([y,\infty))\leq C\,y^{-(\alpha-1)}.
\ee

Let $h$ and $g$ be as in the statement of the proposition, and
let $N$ be an integer such that $2^{-N}\leq a$. Then, using
(\ref{lower-tech5}),
$$\sum_{n=N}^\infty \Theta_\alpha\Big(\int \ell^a(d\sigma)
\,{\bf 1}_{\{\m(B(\sigma,2^{-n}))\geq h(2^{-n})\}}\Big)
=\sum_{n=N}^\infty \mu([g(2^{-n}),\infty))\leq C\sum_{n=N}^\infty
g(2^{-n})^{-(\alpha-1)}<\infty$$
by our assumption (\ref{sum-condi}). Hence,
$$\sum_{n=N}^\infty {\bf 1}_{\{\m(B(\sigma,2^{-n}))\geq
h(2^{-n})\}}<\infty\;,\quad
\ell^a(d\sigma)\hbox{ a.e., }\Theta_\alpha\hbox{ a.e.}$$
and so
$$\limsup_{n\to\infty} {\m(B(\sigma,2^{-n}))\over h(2^{-n})}\leq 1
\;,\quad
\ell^a(d\sigma)\hbox{ a.e., }\Theta_\alpha\hbox{ a.e.}$$
Since this holds for every $a>0$, we can replace $\ell^a(d\sigma)$ a.e.
by $\m(d\sigma)$ a.e. in the last display. By Lemma
\ref{comparison-Hausdorff} (i), this
implies
$$h-m(\t)>0\;,\quad\Theta_\alpha\hbox{ a.e.}$$

Finally, we may find a function $\wt h(r)=r^{\alpha\over \alpha-1}\wt
g(r)$, such that $\wt h$ and $\wt g$ satisfy the same assumptions as
$h$ and $g$, and $\wt g(r)/g(r) \la 0$ as $r\to 0$.
We have $\wt h-m(\t)>0$ which implies $h-m(\t)=\infty$. This completes
the proof of Proposition \ref{lower-Haus}. \cq

\medskip
We now turn to the proof of part (ii) of Theorem \ref{stable-Haus}.
We thus fix $u<0$, and we aim at proving that $h_u-m(\t)=0$,
$\Theta_\alpha$ a.e. We also fix $\delta\in(0,1/2)$ and an integer
$n_0\geq 1$ such that $2^{-n_0}<\delta$. The main step of the
proof is to control the Hausdorff measure $h_u-m(B_{n_0})$
of the ``bad set''
$$B_{n_0}=\{\sigma\in\t: 2\delta\leq d(\rho,\sigma)\leq (2\delta)^{-1}
\hbox{ and } \m(B(\sigma,2^{-n}))\leq h_u(2^{-n})
\hbox{ for every }n\geq n_0\}.$$

Let $p\geq n_0+3$ be an integer. For every integer $k\geq 1$ denote
by
$(\t^{k,p}_j)_{1\leq j\leq N_{k,p}}$ the subtrees of $\t$
above level $k2^{-p}$ with height greater than $2^{-p}$ (cf
Section 4.2 in \cite{DuLG}).
Also set
$$\wt\t^{k,p}_j=\{\sigma\in\t^{k,p}_j: d(\rho,\sigma)\in[k2^{-p},
(k+2)2^{-p})\}.$$

To simplify notation, we put
$$I_p=\{(k,j):k\geq 1,1\leq j\leq N_{k,p}\}.$$
Suppose that $\wt\t^{k,p}_j\cap B_{n_0}\not =\varnothing$ for some
$(k,j)\in I_p$, and let $\sigma_0\in
\wt\t^{k,p}_j\cap B_{n_0}$. Then, for every $\sigma\in\wt\t^{k,p}_j$,
$$B(\sigma,2^{-n})\subset B(\sigma_0,2^{-n}+4\,2^{-p})\subset
B(\sigma_0,2^{-n+1})$$
provided that $n\leq p-2$. Since $\sigma_0\in B_{n_0}$, we have,
for every $\sigma\in \wt\t^{k,p}_j$,
$$\m(B(\sigma,2^{-n}))\leq h_u(2^{-n+1})\;,\hbox{ for every }
n\in\{n_0+1,n_0+2,\ldots,p-2\}.$$
Thus if we set
\ba &&\wt B_{n_0,p}=\{\sigma\in\t: \delta\leq d(\rho,\sigma)\leq
\delta^{-1}
\hbox{ and } \m(B(\sigma,2^{-n}))\leq h_u(2^{-n+1})\\
&&\hspace{7cm}\hbox{ for every }n\in\{n_0+1,n_0+2,\ldots,p-2\}\},
\ea
we see that the condition $\wt\t^{k,p}_j\cap B_{n_0}\not =\varnothing$
implies $\wt\t^{k,p}_j\subset \wt B_{n_0,p}$.

It follows that, for every real $b>0$,
\begin{eqnarray}
\label{upper-tech1}
&&\!\!\!\!\#\{(k,j)\in I_p:\wt\t^{k,p}_j\cap B_{n_0}\not
=\varnothing\}\nonumber\\ &&\leq \#\{(k,j)\in I_p:k\leq 2^p\delta^{-1}
\hbox{ and }\m(\wt\t^{k,p}_j)<b\} +
b^{-1}\sum_{k=1}^\infty\sum_{j=1}^{N_{k,p}}
\int_{\wt\t^{k,p}_j}\m(d\sigma)\,
{\bf 1}_{\wt B_{n_0,p}}(\sigma)\nonumber\\
&&\leq \#\{(k,j)\in I_p:k\leq 2^p\delta^{-1}
\hbox{ and }\m(\wt\t^{k,p}_j)<b\} +2 b^{-1}\int \m(d\sigma)\,{\bf 1}_{\wt
B_{n_0,p}}(\sigma).
\end{eqnarray}
In the last bound we used the fact that, for every $\sigma\in\t$, there
are at most two pairs $(k,j)\in I_p$ such that $\sigma\in \wt\t^{k,p}_j$.

We will apply the bound (\ref{upper-tech1}) with
$$b=b_p=2^{-p\,{\alpha\over \alpha-1}}\,(\log p)^{-\kappa},$$
where $\kappa>{\alpha\over \alpha-1}$ is arbitrary. We use
different arguments to bound the two terms in the
right-hand side of (\ref{upper-tech1}). To bound the first term,
we apply Lemma \ref{technical-stable}. From this lemma and the
scaling properties of the stable tree recalled in the
preceding section, we have, for every $b>0$ and $r>0$,
\be
\label{upper-tech2}
\Theta_\alpha\Big(\m(\t_{\leq r})\leq b\;,\;H(\t)\geq r\Big)
\leq C_\alpha\,r^{-{1\over \alpha-1}}\,\exp(-c_\alpha
r\,b^{-{\alpha-1\over
\alpha}}).
\ee
On the other hand, the branching property of the stable tree
(cf Theorem 4.2 in \cite{DuLG}) guarantees that for every $k\geq 1$,
under $\Theta_\alpha(\cdot\mid H(\t)\geq k2^{-p})$ and
conditionally given $\langle\ell^{k2^{-p}},1\rangle$, the
trees $\t^{k,p}_1,\ldots,\t^{k,p}_{N_{k,p}}$ are distributed as
the atoms of a Poisson point measure with intensity
$\langle\ell^{k2^{-p}},1\rangle\,\Theta_\alpha(\cdot\cap\{H(\t)\geq
2^{-p})$. Recalling that
$\Theta_\alpha(\langle\ell^{k2^{-p}},1\rangle)=1$, we get
\ba
&&\Theta_\alpha\Big(\#\{(k,j)\in I_p:k\leq 2^p\delta^{-1}
\hbox{ and }\m(\wt\t^{k,p}_j)<b_p\}\Big)\\
&&\qquad=\sum_{k=1}^{[\delta^{-1}2^p]}
\Theta_\alpha\Big(\#\{j:1\leq j\leq N_{k,p}
\hbox{ and }\m(\wt\t^{k,p}_j)<b_p\}\Big)\\
&&\qquad=\sum_{k=1}^{[\delta^{-1}2^p]}
\Theta_\alpha(\m(\t_{\leq 2\cdot2^{-p}})\leq b_p,\,H(\t)\geq 2^{-p})
\\
&&\qquad\leq \delta^{-1}2^p\,
\Theta_\alpha(\m(\t_{\leq 2^{-p}})\leq b_p,\,H(\t)\geq 2^{-p})\\
&&\qquad\leq C_\alpha\,\delta^{-1}\,2^{p{\alpha\over \alpha-1}}
\exp(-c_\alpha 2^{-p}b_p^{-{\alpha-1\over
\alpha}})
\ea
using (\ref{upper-tech2}) in the last bound. Recalling our choice of
$b_p$, we deduce from this bound that
$$\sum_{p=n_0+3}^\infty
h_u(2^{-p})\,\Theta_\alpha\Big(\#\{(k,j)\in I_p:k\leq 2^p\delta^{-1}
\hbox{ and }\m(\wt\t^{k,p}_j)<b_p\}\Big)<\infty.$$
It follows that
\be
\label{upper-tech3}
\lim_{p\to \infty}
h_u(2^{-p})\,\#\{(k,j)\in I_p:k\leq 2^p\delta^{-1}
\hbox{ and }\m(\wt\t^{k,p}_j)<b_p\}= 0\;,\ \Theta_\alpha
\hbox{ a.e.}
\ee

We now turn to the second term in the right-hand side of
(\ref{upper-tech1}). From the definition of $\wt B_{n_0,p}$
and the identity $\m=\int_0^\infty da\,\ell^a$, we have
$$\Theta_\alpha\Big(\int \m(d\sigma)\,{\bf 1}_{\wt
B_{n_0,p}}(\sigma)\Big)
=\Theta_\alpha\Big(\int_\delta^{\delta^{-1}}
da\int \ell^a(d\sigma)\,{\bf 1}_
{\{\m(B(\sigma,2^{-n}))\leq h_u(2^{-n+1})\;,\;\forall
n\in\{n_0+1,\ldots,p-2\}}\Big).$$
  For every
$n\geq n_0+1$, set
$$\cc(\sigma,2^{-n})=\{\sigma'\in\t: 2^{-n-2}\leq d(\sigma\wedge
\sigma',\sigma)<2^{-n-1}\hbox{ and }0<d(\sigma\wedge \sigma',\sigma')\leq
2^{-n-1}\}.$$
Clearly, $\cc(\sigma,2^{-n})\subset B(\sigma,2^{-n})$ and so
$$\Theta_\alpha\Big(\int \m(d\sigma)\,{\bf 1}_{\wt
B_{n_0,p}}(\sigma)\Big)
\leq\int_\delta^{\delta^{-1}}
da\,\Theta_\alpha\Big(\int \ell^a(d\sigma)\,{\bf 1}_
{\{\m(\cc(\sigma,2^{-n}))\leq h_u(2^{-n+1})\;,\;\forall
n\in\{n_0+1,\ldots,p-2\}}\Big).$$

Let us fix $a\in[\delta,\delta^{-1}]$. It follows from
Theorem \ref{Palmdec} that under the probability measure
$\Theta_\alpha(d\t)\ell^a(d\sigma)$, the random variables
$$\m(\cc(\sigma,2^{-n}))\;,\ n=n_0+1,n_0+2,\ldots$$
are independent, and furthermore the law of
$\m(\cc(\sigma,2^{-n}))$ is determined by
$$\Theta_\alpha\Big(\int \ell^a(d\sigma)\,
\exp(-\lambda \m(\cc(\sigma,2^{-n})))\Big)
=\exp(-\alpha 2^{-n-2}v^0_{2^{-n-1}}(\lambda)^{\alpha-1})$$
where $v^0_r(\lambda)$ is as previously. From the scaling
property (\ref{scaling-v}), we see that the law $\nu$
of
$$2^{n{\alpha\over \alpha-1}}\,\m(\cc(\sigma,2^{-n}))$$
under $\Theta_\alpha(d\t)\ell^a(d\sigma)$ does not depend
on $a$ nor on $n$ (this is indeed true provided $2^{-n-1}\leq a$,
which holds here since $2^{-n-1}\leq 2^{-n_0-1}\leq \delta\leq a$).
Furthermore, the Laplace transform of $\nu$
is
$$\int \nu(dx)\,e^{-\lambda x}=\exp\Big(-{\alpha\over
4}\,v^0_{1/2}(\lambda)^{\alpha-1}\Big).$$
Since $\lambda^{-1}v^0_{1/2}(\lambda)\ua 1/2$ as $\lambda\da 0$,
we get
$$\int \nu(dx)\,e^{-\lambda x}=1-{\alpha\over
2}2^{-\alpha}\,\lambda^{\alpha-1} + o(\lambda^{\alpha-1})$$
as $\lambda\da 0$. From a standard Tauberian theorem, it
follows that there exists a constant $c_0>0$ such that,
for every $y\geq 1$,
$$\nu([y,\infty))\geq c_0\,y^{1-\alpha}.$$

Using this bound together with the previously mentioned independence,
we get
\ba
&&\Theta_\alpha\Big(\int \ell^a(d\sigma)\,{\bf 1}_
{\{\m(\cc(\sigma,2^{-n}))\leq h_u(2^{-n+1})\;,\;\forall
n\in\{n_0+1,\ldots,p-2\}}\Big)\\
&&\quad=\prod_{n=n_0+1}^{p-2} \Big(1-\nu([2^{n{\alpha\over \alpha-1}}
h_u(2^{-n+1}),\infty))\Big)\\
&&\quad
\leq \prod_{n=n_0+1}^{p-2} \Big(1-c_02^{-n\alpha}h_u(2^{-n+1})^{1-\alpha}
\Big)\\
&&\quad=\prod_{n=n_0+1}^{p-2} \Big(1-c_02^{-\alpha}
((n-1)\log 2)^{-1}(\log(n-1)+\log\log 2)^{(1-\alpha)u}\Big)\\
&&\quad\leq\exp\Big(-\ov c_0
\Big((\log(p-2))^{1-(1-\alpha)u}-(\log n_0)^{1-(1-\alpha)u}\Big)\Big)
\ea
where $\ov c_0$ is a positive constant and the last bound
follows from simple analytic estimates.

By integrating with respect to $a$, we arrive at
$$\Theta_\alpha\Big(\int \m(d\sigma)\,{\bf 1}_{\wt
B_{n_0,p}}(\sigma)\Big)
\leq \delta^{-1}\,\exp\Big(-\ov c_0
\Big((\log(p-2))^{1-(1-\alpha)u}-(\log n_0)^{1-(1-\alpha)u}\Big)\Big).$$
Notice that $1-(1-\alpha)u>1$ since $u<0$. It then follows from the
preceding bound that
$$\sum_{p=n_0+3}^\infty h_u(2^{-p})b_p^{-1}\,
\Theta_\alpha\Big(\int \m(d\sigma)\,{\bf 1}_{\wt
B_{n_0,p}}(\sigma)\Big)<\infty$$
and thus
\be
\label{upper-tech4}
\lim_{p\to\infty} h_u(2^{-p})b_p^{-1}\,
\int \m(d\sigma)\,{\bf 1}_{\wt
B_{n_0,p}}(\sigma)=0\;,\ \Theta_\alpha\hbox{ a.e.}
\ee

By (\ref{upper-tech1}), (\ref{upper-tech3}) and (\ref{upper-tech4}),
we have
$$\lim_{p\to\infty} h_u(2^{-p})
\#\{(k,j)\in I_p:\wt\t^{k,p}_j\cap B_{n_0}\not
=\varnothing\}=0\;,\ \Theta_\alpha\hbox{ a.e.}$$
Since the sets $\wt\t^{k,p}_j\cap B_{n_0}$ provide a covering
of $B_{n_0}$ by sets with diameter less than $4\,2^{-p}$, the definition
of Hausdorff measure gives
$$h_u-m(B_{n_0})=0\;,\ \Theta_\alpha\hbox{ a.e.}$$
By passing to the limit $n_0\ua \infty$ and $\delta\da 0$, we obtain
$$h_u-m\Big(\Big\{\sigma\in \t:\limsup_{n\to \infty} {
\m(B(\sigma,2^{-n}))\over h_u(2^{-n})}<1\Big\}\Big)=0
\;,\ \Theta_\alpha\hbox{ a.e.}$$
On the other hand, Lemma \ref{comparison-Hausdorff} (ii) yields
$$h_u-m\Big(\Big\{\sigma\in \t:\limsup_{n\to \infty} {
\m(B(\sigma,2^{-n}))\over h_u(2^{-n})}\geq 1\Big\}\Big)
\leq M_2\,\m(\t)<\infty
\;,\ \Theta_\alpha\hbox{ a.e.}$$
We conclude that $h_u-m(\t)<\infty$
and since this holds for every $u<0$, we must indeed have
$h_u-m(\t)=0$, $\Theta_\alpha$ a.e. \cq

\medskip
{\bf Proof of Theorem \ref{stable-local-Haus}}. Many arguments here
are similar to the preceding proof, and we will only sketch details.
Without loss of generality we may take $a=1$.
Part (i) is a consequence of the following proposition.

\begin{proposition}
\label{lower-stable-Haus}
Suppose that $\wt h:[0,\infty)\to[0,\infty)$ is a monotone increasing
function that can be written in the form $\wt h(r)=r^{1\over
\alpha-1}\wt g(r)$ where $\wt g$ is monotone decreasing in a neighborhood of
the origin. Then the condition
\be
\label{sum-condi-local}
\sum_{n=1}^\infty \wt g(2^{-n})^{-(\alpha-1)}<\infty
\ee
implies that $\wt h-m(\t(1))=\infty$, $\Theta_\alpha$ a.e. on $\{H(\t)>1\}$.
\end{proposition}

\noindent{\bf Conjecture.} If (\ref{sum-condi-local}) fails, then
$\wt h-m(\t(1))=0$, $\Theta_\alpha$ a.e.

\smallskip
\proof Using Theorem \ref{Palmdec} in the same way as in the
derivation of (\ref{lower-tech3}),
we  have for $\varepsilon\in(0,1]$ and $\lambda>0$,
$$\Theta_\alpha\Big(\int
\ell^1(d\sigma)\,\exp-\lambda\ell^1(B(\sigma,\varepsilon))\Big)
=\exp-\alpha\int_0^{\varepsilon/2}dr\,u_r(\lambda)^{\alpha-1}$$
where $u_r(\lambda)=\Theta(1-\exp-\lambda\langle\ell^r,1\rangle)$ is
given by (\ref{Laplexp}). Straightforward calculations now give
$$\Theta_\alpha\Big(\int
\ell^1(d\sigma)\,\exp-\lambda\ell^1(B(\sigma,\varepsilon))\Big)
=\Big(1+{(\alpha-1)\lambda^{\alpha-1}\varepsilon\over
2}\Big)^{-{\alpha\over \alpha-1}}.$$
Hence the law $\wt\mu$ of $\varepsilon^{-{1\over
\alpha-1}}\ell^1(B(\sigma,\varepsilon))$
under $\Theta_\alpha(d\t)\ell^1(d\sigma)$ does not depend on
$\varepsilon\in(0,1]$ and as
in the proof of Proposition \ref{lower-Haus}, there is a constant $\wt C$ such
that, for every $y>0$,
$$\wt\mu([y,\infty))\leq \wt C\,y^{-(\alpha-1)}.$$
If $\wt h$ is as in the statement of the proposition, it follows that
$$\limsup_{n\to\infty} {\ell^1(B(\sigma,2^{-n}))\over \wt h(2^{-n})}\leq 1
\;,\quad
\ell^1(d\sigma)\hbox{ a.e., }\Theta_\alpha\hbox{ a.e.}$$
The end of the proof is now similar to that of Proposition
\ref{lower-Haus}. \cq

\medskip
Let us now turn to the proof of (ii). The outline is again similar to the proof
of Theorem \ref{stable-Haus} (ii) but there are a few minor differences.
We fix $u<0$ and an integer $n_0\geq 1$. The ``bad set'' is now defined by
$$B_{n_0}=\{\sigma\in\t(1):\ell^1(B(\sigma,2^{-n}))\leq \wt h_u(2^{-n})
\hbox{ for every }n\geq n_0\}.$$
If $p\geq n_0+2$ is an integer, we denote by $\t^p_j$, $1\leq j\leq N_p$
the subtrees of $\t$ above level $1-2^{-p}$ that intersect $\t(1)$.
Arguing in the proof
of Theorem \ref{stable-Haus}, we can check that if $\t^p_j\cap B_{n_0}\not =
\varnothing$, then $\t^p_j\cap\t(1)\subset \wt B_{n_0,p}$, where
$$\wt B_{n_0,p}=\{\sigma\in\t(1):\ell^1(B(\sigma,2^{-n}))\leq h_u(2^{-n+1})
\hbox{ for every }n\in\{n_0+1,\ldots,p-1\}\}.$$
It follows that, for every $b>0$,
\be
\label{stable-tech1}
\#\{j\leq N_p:\t^p_j\cap B_{n_0}\not =\varnothing\}
\leq \#\{j\leq N_p:\ell^1(\t^p_j)<b\} +b^{-1}\int \ell^1(d\sigma)\,{\bf 1}_{\wt
B_{n_0,p}}(\sigma).
\ee
We apply this estimate with $b=b_p=2^{-{p\over \alpha-1}}p^{-\kappa}$, where
$\kappa>\alpha(\alpha-1)^{-2}$.

To bound the first term in the right-hand side of (\ref{stable-tech1}), we use
(\ref{Laplexp}) to get for every $\lambda>0$ and $r>0$,
$$\Theta_\alpha(\exp-\lambda \langle\ell^r,1\rangle \mid
\langle\ell^r,1\rangle >0)
=1-\Big({(\alpha-1)r+\lambda^{1-\alpha}\over
(\alpha-1)r}\Big)^{1\over 1-\alpha}.$$
It follows that there is a constant $C_\alpha$ such that, for every $r>0$
and $b>0$,
$$\Theta_\alpha(0<\langle\ell^r,1\rangle\leq b)\leq C_\alpha
\,r^{-{\alpha\over \alpha-1}}
\,b^{\alpha-1}.$$
Using the branching property as in the proof
of Theorem \ref{stable-Haus}, we get
$$\Theta_\alpha(\#\{j\leq N_p:\ell^1(\t^p_j)<b_p\})
=\Theta_\alpha(0<\langle\ell^{2^{-p}},1\rangle< b_p)\leq
C_\alpha\,2^{p\over \alpha-1}\,
p^{-\kappa(\alpha-1)},$$
and from the choice of $\kappa$, we have
\be
\label{stable-tech6}
\lim_{p\to\infty}\wt h_u(2^{-p})\,\#\{j\leq
N_p:\ell^1(\t^p_j)<b_p\}=0\;,\ \Theta_\alpha
\hbox{ a.e.}
\ee

In order to bound the second term in the right-hand side of
(\ref{stable-tech1}), we set for every $\sigma\in\t(1)$ and every
integer $n\geq 1$,
$$\wt{\cal C}(\sigma,2^{-n})=\{\sigma'\in\t(1):1-2^{-n-1}<
d(\rho,\sigma\wedge \sigma')
\leq 1-2^{-n-2}\}$$
in such a way that $\wt{\cal C}(\sigma,2^{-n})\subset
B(\sigma,2^{-n})$. It easily
follows from Theorem \ref{Palmdec} that the random variables $\ell^1(\wt{\cal
C}(\sigma,2^{-n}))$, $n\geq1$ are independent under the
probability measure $\Theta_\alpha(d\t)\ell^1(d\sigma)$. Furthermore,
simple calculations give for every $\lambda>0$,
\ba
\Theta_\alpha\Big(\int \ell^1(d\sigma)\,\exp-\lambda\ell^1(\wt{\cal
C}(\sigma,2^{-n}))\Big)
&=&\exp\Big(-\alpha\int_{2^{-n-2}}^{2^{-n-1}}dr\,u_r(\lambda)^{\alpha-1}\Big)\\
&=&\Big({\lambda^{1-\alpha}+(\alpha-1)2^{-n-1}\over\lambda^{1-\alpha}+
(\alpha-1)2^{-n-2}}\Big)
^{-{\alpha\over \alpha-1}}.
\ea
Hence the law $\wt\nu$ of $2^{n\over\alpha-1}\ell^1(\wt{\cal
C}(\sigma,2^{-n}))$ under $\Theta_\alpha(d\t)\ell^1(d\sigma)$
does not depend on $n$. From the preceding Laplace transform, we also get
the existence of a constant $\wt c_0>0$ such that, for every $y\geq 1$,
$$\wt\nu([y,\infty))\geq \wt c_0\,y^{1-\alpha}.$$
Using this lower bound and the previously mentioned independence, the
same calculations as
in the proof of Theorem \ref{stable-Haus} lead to
$$\Theta_\alpha\Big(\int \ell^1(d\sigma)\,{\bf 1}_{\wt B_{n_0,p}}(\sigma)\Big)
\leq \exp\Big(-c'_0\Big((\log(p-1))^{1-(1-\alpha)u}-(\log
n_0)^{1-(1-\alpha)u}\Big)\Big)$$
where $c'_0$ is a positive constant. Since $1-(1-\alpha)u>1$, it
easily follows that
\be
\label{stable-tech7}
\lim_{p\to\infty} \wt h_u(2^{-p})b_p^{-1}\int \ell^1(d\sigma)\,{\bf
1}_{\wt B_{n_0,p}}(\sigma)
=0\;,\ \Theta_\alpha\hbox{ a.e.}
\ee
Thanks to (\ref{stable-tech1}), (\ref{stable-tech6}) and
(\ref{stable-tech7}), the remaining
part
of the proof is now similar to the end of the proof of Theorem
\ref{stable-Haus}. \cq

\end{document}